\documentclass[a4paper,12pt]{amsart}

\usepackage{latexsym}
\usepackage{amscd}
\usepackage{graphicx}
\usepackage{a4wide}
\usepackage{amssymb}
\usepackage{euscript}
\usepackage{subfigure}
\usepackage{verbatim}
\usepackage{varioref}
\usepackage{wrapfig}
\usepackage[rflt]{floatflt}
\usepackage{epsfig}
\usepackage{rotating}
\usepackage{fancyhdr}

\title{Boundary slopes and the logarithmic limit set}

\author[Stephan Tillmann]{Stephan Tillmann}
\date{\today}

\def\bkC{{\rm \kern.24em \vrule width.05em height1.4ex depth-.05ex \kern-.26em C}}
\def\C{\bkC}
\def\bksC{{\rm \kern.24em \vrule width.05em height1ex depth-.05ex \kern-.26em C}}
\def\sC{\bksC}
\def\bkE{{\rm I\kern-.22em E}}
	
\def\bkH{{\rm I\kern-.22em H}}
			
\def\bkQ{{\rm \kern.24em \vrule width.05em height1.4ex depth-.05ex \kern-.26em Q}}

\def\bkR{{\rm I\kern-.17em R}}
\def\RR{\bkR}
\def\bkZ{{\rm Z\kern-.32em Z}}
\def\Z{\bkZ}
\def\bksZ{{\rm Z\kern-.22em Z}}
\def\sZ{\bksZ}

\def\SL{SL_2(\C)}
\def\PSL{PSL_2(\C)}

\def\Ei{\mathfrak{E}}
\def\PEi{\overline{\Ei}}

\def\k{\mathfrak{k}}

\def\X{\mathfrak{X}}
\def\PX{\overline{\mathfrak{X}}}
\def\R{\mathfrak{R}}
\def\PR{\overline{\mathfrak{R}}}
\def\Red{\mathfrak{Red}}

\def\BC{\mathfrak{BC}}

\def\P{\mathcal{P}}

\def\tree{\mathcal{T}}

\def\l{\mathcal{L}}
\def\m{\mathcal{M}}

\def\prho{\overline{\rho}}

\def\G{\Gamma}

\def\whl{\mathcal{W}} 

\DeclareMathOperator{\Hom}{Hom}
\DeclareMathOperator{\tr}{tr}

\DeclareMathOperator{\im}{im}

\DeclareMathOperator{\te}{t}
\DeclareMathOperator{\pte}{\overline{t}}

\DeclareMathOperator{\qe}{q}

\theoremstyle{plain}

\newtheorem{thm}{Theorem}
\newtheorem*{thm*}{Theorem}
\newtheorem{lem}[thm]{Lemma}
\newtheorem*{lem*}{Lemma}
\newtheorem{cor}[thm]{Corollary}
\newtheorem*{cor*}{Corollary}

\newtheorem*{cla*}{Claim}
\newtheorem{pro}[thm]{Proposition}
\newtheorem*{pro*}{Proposition}

\newtheorem*{fac*}{Fact}

\newtheorem*{que*}{Question}

\newtheorem*{con*}{Conjecture}

\newtheorem*{rem*}{Remarks}
\newtheorem{defn}[thm]{Definition}
\newtheorem*{defn*}{Definition}

\addtocounter{MaxMatrixCols}{4}

\begin{document}

\begin{abstract}
The $A$--polynomial of a manifold whose boundary consists of a single torus
is generalised to an \emph{eigenvalue variety} of a manifold whose boundary
consists of a finite number of tori, and the set of strongly detected
boundary curves is determined by Bergman's 
\emph{logarithmic limit set}, which describes the exponential behaviour of
the eigenvalue variety at infinity.
This enables one to read off the detected boundary curves of a multi--cusped
manifold in a similar way to the 1--cusped case, where the slopes are
encoded in the Newton polygon of the $A$--polynomial. 
\end{abstract}

\maketitle


\section*{Introduction}

Let $V$ be an algebraic variety in $(\C-\{0\})^m$. The first section
discusses the \emph{logarithmic limit set} $V_\infty$ of $V$, which we
use to define the set of ideal points of $V$. We show that each point
with rational coordinate ratios in $V_\infty$ is an ideal point of a curve in
$V$. Thus, the logarithmic limit set lends itself to applications of
Culler--Shalen theory, where essential surfaces are associated to ideal
points of curves in the character variety of a 3--manifold.

Given a compact, orientable, irreducible 3--manifold with boundary
consisting of a single torus, one calls boundary slopes of essential surfaces
associated to ideal points of the character variety \emph{strongly
detected}. It is shown in \cite{ccgls} that the slope of each side of the
Newton polygon of the $A$--polynomial is a strongly detected boundary slope.

In the second section, we generalise the $A$--polynomial to define an
\emph{eigenvalue variety} of a 
manifold with boundary consisting of $h$ tori. We describe the \emph{strongly
detected boundary curves}, i.e. the boundary curves of essential surfaces
associated to ideal points of the character variety, via a map from the
logarithmic limit set of the eigenvalue variety to the
$(2h-1)$--sphere.

The third section contains a computation of the $A$--polynomials of torus
knots. These seem to be the only infinite families of $A$--polynomials
currently known.


\section[Logarithmic limit set]{Bergman's logarithmic limit set}
\label{section:log lim set}


\subsection{Logarithmic limit set}
\label{combinatorics:Logarithmic limit set}

Let $V$ be a subvariety of $(\C -\{ 0\})^m$ defined by an ideal $J$, and
let $\C [X^{\pm}]=\C [X_1^{\pm 1}, \ldots , X_m^{\pm 1}]$. Bergman gives the
following three descriptions of a \emph{logarithmic limit set} in \cite{be}:

1. Define $V_\infty^{(a)}$
as the set of limit points on $S^{m-1}$ of the set of
real $m$--tuples in the interior of $B^m$:
\begin{equation} \label{log lim def1}
  \bigg\{ \frac{ (\log |x_1|, \ldots , \log |x_m|)}
            {\sqrt{1 + \sum (\log |x_i|)^2}}
	              \mid x \in V \bigg\} .
\end{equation}

2. Define $V_\infty^{(b)}$ as the set of $m$--tuples
\begin{equation} \label{log lim def2}
   (-v(X_1), \ldots , -v(X_m))
\end{equation}
as $v$ runs over all real--valued valuations on $\C [X^{\pm}]/J$ satisfying
$\sum v(X_i)^2 = 1$.

3. Define the support $s(f)$ of an element $f \in \C [X^{\pm}]$ to be the
set of points $\alpha = (\alpha_1,{\ldots} ,\alpha_m) \in \Z^m$ such that
$X^\alpha = X_1^{\alpha_1} \cdots X_m^{\alpha_m}$ occurs with non--zero
coefficient in $f$. Then define $V_\infty^{(c)}$ to be the set of 
$\xi \in S^{m-1}$ such that for all non--zero $f \in J$, the maximum value
of the dot product $\xi \cdot \alpha$ as $\alpha$ runs over $s(f)$ is
assumed at least twice.

\begin{figure}[h]
    \includegraphics[width=5cm]{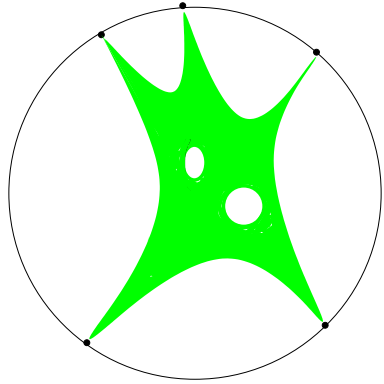}
  \caption{A logarithmic limit set in $S^1$}
  \label{fig:loglim}
\end{figure}

Bergman shows in \cite{be} that the second and third descriptions are
equivalent, and it follows from Bieri and Groves \cite{bg} that all
descriptions are equivalent:
\begin{thm}[Bergman, Bieri--Groves]
$V_\infty^{(a)}=V_\infty^{(b)}=V_\infty^{(c)}$.
\end{thm}
We therefore write $V_\infty$ for the logarithmic limit set of the variety
$V$. A generic picture for the projection of a variety in
$(\C - \{ 0\})^2$ into $B^2$ and its logarithmic limit set in $S^1$ is given
in Figure \ref{fig:loglim}. 
The ideal generated by a set of polynomials $\{ f_i\}_i$ shall be
denoted by $I(f_i)_i$, and the variety generated by the ideal $I$ is denoted
by $V(I)$ or $V(f_i)_i$, if $I = I(f_i)_i$.


\subsection{Newton Polytopes}

The third description of the logarithmic limit set can be
illustrated using spherical duals of Newton polytopes.
Let $f \in \C[X^\pm ]$ be given as an expression
$f = \sum_{\alpha \in \sZ^m} a_\alpha X^\alpha$. Then the
\emph{Newton polytope} of $f$ is the convex hull of
$s(f)=\{ \alpha \in \Z^m \mid a_\alpha \ne 0\}$. Thus, if
$s(f) = \{ \alpha_1,{\ldots} ,\alpha_k\}$, then
\begin{equation}
Newt(f) = Conv(s(f)) = 
\bigg\{ \sum_{i=1}^{k} \lambda_i \alpha_i \mid \lambda \ge
0, \sum_{i=1}^{k} \lambda_i =1 \bigg\}.
\end{equation}
Note that $Newt(0) = \emptyset$, and since
$V(0) = (\C-\{ 0 \} )^m$, we have $V(0)_\infty = S^{m-1}$.

For arbitrary subsets $A,B$ of $\RR^m$, let
$A+B = \{ a+b \mid a \in A, b \in B\}$. The following facts about Newton
polytopes of polynomials are well known (cf. \cite{clo}):

\begin{enumerate}
\item Let $f,g \in \C[X^\pm ]$ such that $fg \ne 0$. Then
$Newt(fg) = Newt(f) + Newt(g)$.
\item $Newt(f+g) \subseteq Conv( s(f) \cup s(g) ) = Conv (Newt(f) \cup
              Newt(g))$. 
\end{enumerate}


\subsection{Spherical Duals}
\label{Spherical Duals}

The \emph{spherical dual} of a bounded convex polytope $P$ in $\RR^m$ is the
set of vectors $\xi$ of length $1$ such that the supremum
$\sup_{\alpha \in P}\alpha \cdot \xi$ is achieved for more than one
$\alpha$, and it is denoted by $Sph(P)$. Geometrically, $Sph(f)$ consists of
all outward pointing unit normal vectors to the support planes of $P$ which
meet $P$ in more than one point. If $P$ is the Newton polytope of a
non--zero polynomial $f$, then the spherical dual of $Newt(f)$ is also
denoted by $Sph(f)$. The following lemma is immediate from the third
description of the logarithmic limit set:

\begin{lem}\label{log:int lem}
Let $V$ be an algebraic variety in $(\C-\{ 0\})^m$ defined by the ideal $J$.
Then $V_\infty$ is the intersection over all non--zero elements of $J$ of the
spherical duals of their Newton polytopes.
\end{lem}

Spherical duals are easy to visualise in low dimensions. Given a convex
polygon $P$ in $\RR^2$, its spherical dual is the collection of points on
the unit circle defined by outward pointing unit normal vectors to edges of
$P$. Given a convex polyhedron in $\RR^3$, we obtain the vertices of its
spherical dual again as points on the unit sphere $S^2$ arising from outward
pointing unit normal vectors to faces. We then join two of these points
along the shorter geodesic arc in $S^2$ if the corresponding faces have a
common edge. This gives a finite graph in $S^2$.


\subsection{}

If the ideal $J$ is principal, defined by an element $f \ne 0$, then Bergman
states in \cite{be} that the set $V(f)_\infty$ is precisely $Sph(f)$.
We give a proof using the following lemma.

\begin{lem}\label{log:sph lem}
Let $f,g \in \C[X^\pm ]$. If $fg \ne 0$, then $Sph(fg) = Sph(f) \cup Sph(g)$.
\end{lem}

\begin{proof}
If $\xi \in Sph(f)$, then there are distinct 
$\alpha_0, \alpha_1 \in Newt(f)$ such that 
\begin{equation*}
\alpha_0 \cdot \xi = \alpha_1 \cdot \xi = \sup_{\alpha \in Newt(f)}
\alpha\cdot \xi.
\end{equation*}
Since $Newt(g)$ is convex and bounded, there is $\beta_0 \in Newt(g)$ such
that 
\begin{equation*}
\beta_0 \cdot \xi = \sup_{\beta \in Newt(g)} \beta\cdot \xi.
\end{equation*}
Since $Newt(fg) = Newt(f) + Newt(g)$, it follows that
$\beta_0 + \alpha_0, \beta_0 + \alpha_1 \in Newt(fg)$ and
\begin{equation*}
(\alpha_0 + \beta_0)\cdot \xi = (\alpha_1+\beta_0) \cdot \xi = 
\sup_{\gamma \in Newt(fg)} \gamma\cdot \xi.
\end{equation*}
Thus, $Sph(f) \subseteq Sph(fg)$, and by symmetry $Sph(g) \subseteq Sph(fg)$.

Now then assume that $\xi \in Sph(fg)$. Then there are distinct $\gamma_0$
and $\gamma_1$ in $Newt(fg)$ such that 
\begin{equation*}
\gamma_0 \cdot \xi = \gamma_1 \cdot \xi = \sup_{\gamma \in Newt(fg)}
\gamma\cdot \xi.
\end{equation*}
Furthermore, there are $\alpha_i \in Newt(f)$ and $\beta_i \in Newt(g)$ such
that $\gamma_i = \alpha_i + \beta_i$. If
$\alpha_0 \cdot \xi  < \alpha_1 \cdot \xi$, then
$\gamma_0 \cdot \xi < (\alpha_1 + \beta_0) \cdot \xi$
and $\alpha_1 + \beta_0 \in Newt(fg)$, contradicting the
choice of $\gamma_0$. Thus, $\alpha_0 \cdot \xi = \alpha_1 \cdot \xi$, and
similarly $\beta_0 \cdot \xi = \beta_1 \cdot \xi$. A similar argument shows:
\begin{equation*}
\alpha_0 \cdot \xi = \alpha_1 \cdot \xi = \sup_{\alpha \in Newt(f)}
\alpha\cdot \xi,
\end{equation*}
and similarly for the $\beta_i$. Since $\gamma_0 \ne \gamma_1$, either
$\alpha_0 \ne \alpha_1$ or $\beta_0 \ne \beta_1$, and hence
either $\xi \in Sph(f)$ or $\xi \in Sph(g)$.
\end{proof}

\begin{pro}
Let $f \in \C[X^\pm ]$. If $f \ne 0$, then $V(f)_\infty = Sph(f)$.
\end{pro}

\begin{proof}
Since each non--zero element of $I(f)$ is of the form $fg$ for some non--zero
$g \in \C[X^\pm ]$, Lemma \ref{log:int lem} and Lemma \ref{log:sph lem} yield:
\begin{equation*}
V(f)_\infty = \underset{0 \ne g \in \sC[X^\pm ]}{\bigcap} Sph(fg)
            = \underset{0 \ne g \in \sC[X^\pm ]}{\bigcap} (Sph(f) \cup Sph(g))
	    = Sph(f),
\end{equation*}
since $Sph(1) = \emptyset$.
\end{proof}

Logarithmic limit sets of varieties defined by a principal ideal are
therefore completely determined by a generator of the ideal, and readily
computable.
If the ideal $J$ is not principal, then $V(J)_\infty$ is the intersection of
the spherical duals of the Newton polytopes of \emph{finitely} many elements
in $J$. This follows from Bieri and Groves \cite{bg}.
However, there is presently no known algorithm to determine a suitable finite
set of elements of $J$. Conjecturally, a Gr\"obner basis may suffice.


\subsection{}
\label{log:computing}

As Bergman notes, the spherical dual of the convex hull of a finite subset
$A\subset \Z^m$ of cardinality $r$ is a finite union of convex
spherical polytopes. It is the union over all $\alpha_0, \alpha_1 \in A$ of
the set of $\xi$ satisfying the $2r$ inequalities 
\begin{equation*}
  \alpha_0 \cdot \xi \ge \alpha \cdot \xi 
  \qquad\text{and}\qquad
  \alpha_1 \cdot \xi \ge \alpha \cdot \xi,
\end{equation*}
where $\alpha$ ranges over $A$. This may be useful for calculations, but it
also shows that corners of the convex spherical polytopes which we consider
have rational coordinate ratios, and points with rational coordinate ratios
are dense in the polytope. Such a polytope is called a \emph{rational convex
spherical polytope}.
We summarise further results by Bergman and Bieri--Groves as follows:

\begin{thm}[Bergman, Bieri--Groves] \label{Bergman-Bieri-Groves}
  The logarithmic limit set $V_\infty$ of an algebraic variety $V$ in
  $(\C -\{ 0\})^m$ is a finite union of rational convex spherical polytopes. 
  The maximal dimension for a polytope in this union is achieved by at least
  one member $P$ in this union, and we have
  \begin{equation*}
    \dim_\RR V_\infty = \dim_\RR P = (\dim_\sC V) - 1.
  \end{equation*}
\end{thm}


\subsection{Whitehead link}

\begin{figure}[t]
  \begin{center}
      \includegraphics[width=5cm]{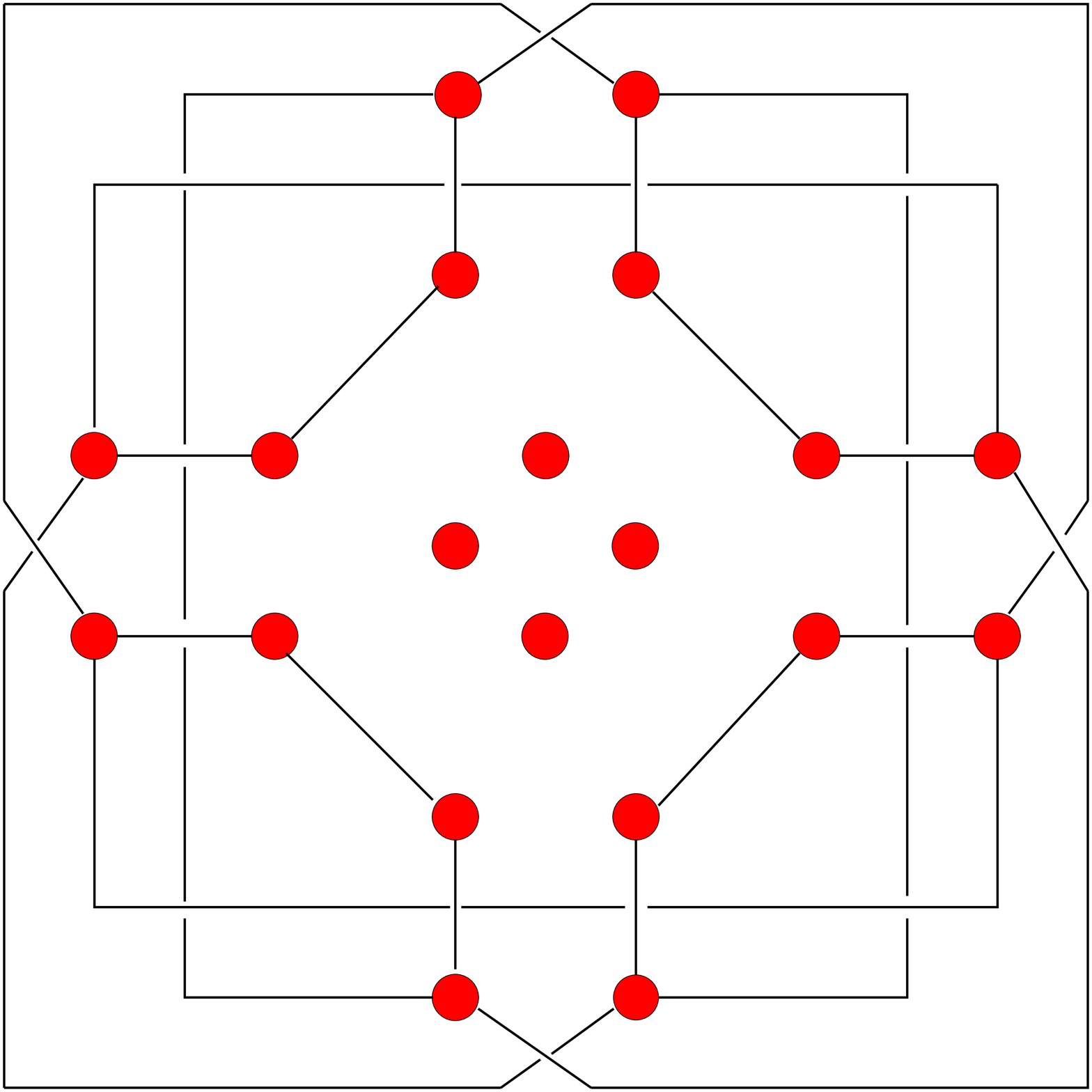}
  \end{center}
  \caption{The logarithmic limit set of $\PEi_0 (\whl)$}
  \label{fig:whl_log_lim}
\end{figure}

Bergman asks in \cite{be} whether every polytope in the above union always
has the same dimension. In \cite{tillus_phd} it is shown that the
irreducible component in the $\PSL$--eigenvalue variety of the Whitehead link
complement $\whl$ which arises from the Dehn surgery component has logarithmic
limit set consisting of five components, four of which are points, and one
of which is a union of geodesic arcs. A schematic picture is given in Figure
\ref{fig:whl_log_lim}.

It would be interesting to
modify Bergman's question and ask whether 
each connected component of the logarithmic limit set of an irreducible
variety can always be expressed as a union of polytopes all having the same
dimension when $\dim_\sC V > 2$.

\subsection{Curve finding lemma}

We conclude this section with a lemma, which is useful for applications of
Culler--Shalen theory. 

\begin{lem}\label{curve finding lemma}
  Let $V$ be an algebraic variety in $(\C -\{ 0\})^m$, and 
  $\xi \in V_\infty$ be a point with rational coordinate ratios. Then
  there is a curve $C$, i.e. an irreducible subvariety of complex
  dimension one, in $V$ such that $\xi \in C_\infty$.
\end{lem}

\begin{proof}
By a lemma of Bergman (\cite{be}, p.464), the point $\xi$ is contained in
the logarithmic limit set of an irreducible component of $V$. Thus, we may
assume that $V$ is irreducible for the purpose of this proof.

Assume that $V$ is defined by the ideal $I$ in 
$\C[X^{\pm}] = \C[X_1^{\pm 1},\ldots, X_m^{\pm 1}]$. We follow word by word Bergman's
construction starting after ``(4)'' on page 465, up to the sentence
``Hence the points $u_y \log |y| \in B^m$ will approach $(0,\ldots,0,-1)$.'' 
on page 466.  Here $u_y$ denotes $(\sqrt{1 + \sum(\log (|y_i|)^2})^{-1}$. Our
adaption reads as follows: 
{\small
\begin{quote}
Let $\xi$ be a point of $V_\infty$ the ratios of whose coordinate functions
are rational. Making use of the action of $GL_m(\Z )$ on $\C[X^{\pm}]$ and
the induced action on $S^{m-1}$ (not isometric!), we can reduce to the case
where $\xi = (0,\ldots ,0,-1)$.

Let $R'$ designate the ring
$\C[X_1^{\pm 1}, \ldots , X_{m-1}^{\pm 1}; X_m] \subset R = \C[X^{\pm}]$, and
$R_0$ the ring $\C[X_1^{\pm 1}, \ldots , X_{m-1}^{\pm 1}]$, into which we map
$R'$ by sending $X_m$ to $0$. Let $I'$ designate $I \cap R'$, and $I_0$ the
image of $I'$ in $R_0$. Thus, $I_0$ consists of the components of degree
$0$ in $X_m$ of elements of $I'$. Let $V'$ designate the subvariety of
$(\C - \{ 0\})^{m-1} \times \C$ defined by $I'$, and
$V_0 = V' \cap (\C - \{ 0\})^{m-1} \times \{ 0\}$, the subvariety of
$(\C - \{ 0\})^{m-1} \times \{ 0\}$ defined by $I_0$. $V'$ will be the
Zariski closure of $V$ in $(\C - \{ 0\})^{m-1} \times \C$, and
$V = V' - V_0$.

If $(0,\ldots ,0,-1) \in V_\infty$, this means that every nonzero $a \in I$
has more than one term of minimal degree in $X_m$. Hence every nonzero
element of $I_0$ has more than one term, i.e. $I_0$ is a proper ideal of
$R_0$. Then $V_0$ is non--empty: let $x \in V_0$.

Now $V$ will in fact be dense in $V'$ under the topology induced by the
absolute value on $\C$. [\ldots ] Hence we can obtain $x$ as a limit of
points $y \in V$. As $y$ approaches $x$, $\log |y_i|$ will approach the
finite values $\log |x_i|$ for $i = 1,\ldots ,m-1$, but
$\log |y_m|$ will approach $-\infty$. Hence the points
$u_y \log |y| \in B^m$ will approach $(0,\ldots ,0,-1)$.
\end{quote}
}

We construct a curve $C$ in $V$ such that $\xi=(0,\ldots ,0,-1)\in C_\infty$
using the following two facts from \cite{mu}: 
\begin{enumerate}
\item (1.14) Let $Y$ be an irreducible affine variety. If $X$ is a proper
subvariety of $Y$, then $\dim X < \dim Y$.
\item (3.14) If $X$ is an r--dimensional affine variety and
$f_1,{\ldots} ,f_k$ are polynomial functions on $X$, then every component of
$X \cap V(f_1,{\ldots} ,f_k)$ has dimension $\ge r-k$.
\end{enumerate}

Since $V_0$ is a proper subvariety of $V'$, obtained by intersection
with a hyperplane, we have $\dim V_0 = \dim V' -1$, and since $V=V'-V_0$, we
have $\dim V = \dim V'$. We need to find a curve $C'$ in $V'$ which has
non--trivial intersection with both $V_0$ and $V$, since then $C=C'\cap V$
is a curve in $V$ such that $\xi\in C_\infty$. 

Let $\dim V = d$. Let $x \in V_0$ and choose a polynomial $f_1$ on $V_0$
such that $f_1(x)=0$, but $f_1$ is not identical to $0$ on $V_0$. Then $f_1$
is not identical to $0$ on $V'$ and each component of the variety 
$V' \cap V(f_1)$ has dimension exactly $d-1$, and intersects $V_0$ in a
$(d-2)$--dimensional subvariety.
Let $V_1$ be a component of $V' \cap V(f_1)$ containing $x$. 

Choose a polynomial $f_2$ such that  $f_2(x)=0$, but $f_2$ is not identical
to $0$ on $V_1 \cap V_0$. Then each component of 
$V_1 \cap V(f_2) = V' \cap V(f_1, f_2)$ has dimension exactly $d-2$ and
intersects $V_0$ in a $(d-3)$--dimensional subvariety. Let $V_2$ be a
component of $V' \cap V(f_1, f_2)$ containing $x$.

We now proceed inductively for $k = 2,{\ldots} ,d-2$. Given a component
$V_k$ of $V' \cap V(f_1,{\ldots} ,f_k)$ which is a $(d-k)$--dimensional
subvariety of $V'$, with the property that it contains $x$, but also points
of $V' -V_0$, we let $f_{k+1}$ be a polynomial function on $V_k$ which
vanishes on $x$ but not on all of $V_k \cap V_0$. Then each component of
$V' \cap V(f_1,{\ldots} ,f_{k+1})$ has dimension exactly $d-k-1$ and
intersects $V_0$ in a $(d-k-2)$--dimensional subvariety, and we let
$V_{k+1}$ be a component of $V' \cap V(f_1,{\ldots} ,f_{k+1})$ containing
$x$.

Then $V_{d-1}\cap V$ is a curve in $V$ with the desired property.
\end{proof}


\section[Eigenvalue variety]{The eigenvalue variety}
\label{Eigenvalue variety}


\subsection{Essential surfaces}
\label{Essential surfaces}

A \emph{surface} $S$ in a compact 3--manifold $M$ will always mean a
2--dimensional PL submanifold \emph{properly embedded} in $M$, that is, a
closed subset of $M$ with $\partial S = S \cap \partial M$. If $M$ is not
compact, we replace it by a compact core. An embedded sphere $S^2$ in a
3--manifold $M$ is called \emph{incompressible} if it does not bound an
embedded ball in $M$, and a manifold is \emph{irreducible} if it contains no
incompressible 2--spheres. 

An orientable surface $S$ without 2--sphere or disc components in an
orientable 3--manifold $M$ is called \emph{irreducible} if for each disc
$D\subset M$ with $D \cap S = \partial D$ there is a disc $D' \subset S$
with $\partial D' = \partial D$. We will also use the following definition:
\begin{defn} \cite{sh1}\label{def:essential}
  A surface $S$ in a compact, irreducible, orientable 3--manifold is said to
    be essential if it has the following five properties:
    \begin{enumerate}
       \item $S$ is bicollared;
       \item the inclusion homomorphism $\pi_1(S_i) \to \pi_1(M)$ is
          injective for
           every component $S_i$ of $S$;
       \item no component of $S$ is a 2--sphere;
       \item no component of $S$ is boundary parallel;
       \item $S$ is nonempty.
    \end{enumerate}
\end{defn}


\subsection{$\mathbf{\partial}$--incompressible}

The following definition and facts can be found in \cite{ha1}. A (properly
embedded) surface $S$ in a compact 3--manifold $M$ is
\emph{$\partial$--incompressible} if for each disc $D\subset M$ such that
$\partial D$ splits into two arcs $\alpha$ and $\beta$ meeting only at their
common endpoints with $D \cap S = \alpha$ and $D \cap \partial M = \beta$
there is a disc $D' \subset S$ with $\alpha \subset \partial D'$ and
$\partial D' - \alpha \subset \partial S$.
A surface $S$ is incompressible and $\partial$--incompressible if and only
if each component of $S$ is incompressible and $\partial$--incompressible. 

\begin{lem} \cite{ha1} \label{ev:bdry incomp}
Let $S$ be a connected incompressible surface in the irreducible 3--manifold
$M$, with $\partial S$ contained in torus boundary components of $M$. 

Then either $S$ is $\partial$--incompressible or $S$ is a boundary parallel
annulus. 
\end{lem}

For the remainder of this section, we let $M$ be a compact, orientable,
irreducible 3--manifold with non--empty boundary consisting of a disjoint
union of $h$ tori. It follows from Definition \ref{def:essential} and
Lemma \ref{ev:bdry incomp} that an essential surface in $M$ is
$\partial$--incompressible.


\subsection{Boundary curve space}
\label{boundary_stuff:Boundary curve space}

Denote the boundary tori of $M$ by $T_1, {\ldots} ,T_h$, and choose an
oriented meridian $\m_i$ and an oriented longitude $\l_i$ for each $T_i$.
Let $S$ be an incompressible and $\partial$-incompressible surface with
non--empty boundary in $M$.
There are coprime integers $p_i, q_i$ such that $S$ meets $T_i$ in
$n_i$ curves parallel to $p_i \m_i + q_i \l_i$. We thus obtain a point  
\begin{equation*}
  (n_1p_1, n_1q_1, \ldots , n_hp_h, n_hq_h) \in \Z^{2h} - \{ 0\}.
\end{equation*}
We view the above point as an element in $\RR P^{2h-1}$ since most methods
of detecting essential surfaces only do so up to projectivisation.
Furthermore, we ignore the orientation of $\partial S$, so the loop
$p_i \m_i + q_i \l_i$ is equivalent to $-p_i \m_i - q_i \l_i$. All we have
to know are the relative signs of $p_i$ and $q_i$. This equivalence induces
a $\Z_2^{h-1}$ action on $\RR P^{2h-1}$, and the closure of the set of
\emph{projectivised boundary curve coordinates}:
\begin{equation*}
    [n_1p_1, n_1q_1, \ldots , n_hp_h, n_hq_h]
      \in \RR P^{2h-1}/\Z_2^{h-1}
\end{equation*}
arising from incompressible and $\partial$--incompressible surfaces with
non--empty boundary in $M$ is called the \emph{boundary curve space} of $M$
and denoted by $\BC (M)$. Note that $\RR P^{2h-1}/\Z_2^{h-1} \cong S^{2h-1}$
(see \cite{ha}). We will use the following result.

\begin{thm} \cite{ha} \label{ev:dim bc}
Let $M$ be a compact, orientable, irreducible 3--manifold with non--empty
boundary consisting of a disjoint union of $h$ tori. 
Then $\dim_\RR \BC (M) < h$.
\end{thm}


\subsection{Character variety}

Let $\G$ be a finitely generated group. The set of representations of
$\G$ into $\SL$ is $\R (\G) = \Hom (\G, \SL )$. This set is called the
\emph{representation variety} of $\G$, and regarded as an affine algebraic
set. Each $\gamma \in \G$ defines a \emph{trace function}
$I_\gamma : \R (\G) \to \C$ by $I_\gamma(\rho) = \tr\rho(\gamma)$, which 
is an element of the coordinate ring $\C [\R (\G)]$.

Two representations are \emph{equivalent} if they differ by
an inner automorphism of $\SL$. A representation is \emph{irreducible} if
the only  subspaces of $\C^2$ invariant under its image are trivial. This is
equivalent to saying that the representation cannot be conjugated to a
representation by upper triangular matrices. Otherwise
a representation is \emph{reducible}.

For each $\rho \in \R (\G)$, its \emph{character} is the
function $\chi_\rho : \Gamma \to \C$ defined by
$\chi_\rho (\gamma ) = \tr\rho (\gamma )$. It turns out that
irreducible representations are determined by characters up to equivalence,
and that the reducible representations form a closed subset $\Red (\G)$ of 
$\R (\G)$.

The collection of characters $\X (\G)$ can be regarded as an affine algebraic
set, which is called the \emph{character variety}. There is a regular map 
$\te : \R (\G) \to \X (\G)$ taking representations to characters. 

If $\G$ is the fundamental group of a topological space $M$, we also write
$\R (M)$ and $\X (M)$ instead of $\R (\G)$ and $\X (\G)$ respectively.

There also is a character variety arising from representations
into $\PSL$, see Boyer and Zhang \cite{bozh}, and the relevant objects are
denoted by placing a bar over the previous notation. The natural map 
$\qe: \X (\G) \to \PX (\G)$ is finite--to--one, but in general not onto. 
As with the $\SL$--character variety, there is a surjective
``quotient'' map $\pte : \PR (\G) \to \PX (\G)$, which is
constant on conjugacy classes, and with the property that if $\prho$ is an
irreducible 
representation, then $\pte^{-1}(\pte (\prho ))$ is the orbit of $\prho$
under conjugation. Again, a representation is irreducible if it is not
conjugate to a representation by upper triangular matrices.


\subsection{Dehn surgery component}
\label{Dehn surgery component}

If a 3--manifold $M$ admits a complete hyperbolic structure of finite
volume, then there 
is a discrete and faithful representation $\pi_1(M) \to \PSL$. This
representation  is neccessarily irreducible, as hyperbolic geometry
otherwise implies that $M$ has infinite volume.
If $M$ is not compact, then a \emph{compact core} of $M$ is a compact
manifold $\overline{M}$ such that $M$ is homeomorphic to the interior of
$\overline{M}$. We will rely on the following result: 

\begin{thm}[Thurston] \cite{t, sh1} \label{thurston thm}
Let $M$ be a complete hyperbolic 3--manifold of finite volume with $h$ cusps,
and let $\prho_0 : \pi_1(M) \to \PSL$ be a discrete and faithful
representation  associated to the complete hyperbolic structure. Then
${\prho_0}$ admits a lift $\rho_0$ into $\SL$ which is still
discrete and faithful. The (unique) irreducible component $\X_0$ in the
$\SL$--character variety containing the character $\chi_0$ of
$\rho_0$ has (complex) dimension $h$. 

Furthermore, if $T_1, \ldots , T_h$ are the boundary tori of a compact core
of $M$, and if $\gamma_i$ is a non--trivial element in $\pi_1(M)$ which is
carried by $T_i$, then
$\chi_0 ( \gamma_i) = \pm 2$ and $\chi_0$ is an isolated point of the set
\begin{equation*}
  X^* = \{ \chi \in \X_0 \mid I^2_{\gamma_1} = \ldots = I^2_{\gamma_h} = 4 \}.
\end{equation*}
\end{thm}

The respective irreducible components containing the so--called
\emph{complete representations} $\rho_0$
and ${\prho_0}$ are denoted by $\R_0(M)$ and $\PR_0(M)$ respectively.
In particular, $\te (\R_0) = \X_0$ and $\pte (\PR_0) = \PX_0$ are called the
respective \emph{Dehn surgery components} of the character varieties of $M$,
since the holonomy 
representations of hyperbolic manifolds or orbifolds obtained by performing
high order Dehn surgeries on $M$ are near $\prho_0$ (see \cite{t}).


\subsection{Eigenvalue variety}
\label{boundary:Eigenvalue variety}

Given a presentation of $\pi_1(M)$ with a finite number, $n$, of
generators, $\gamma_1,{\ldots} ,\gamma_n$, we introduce four affine
coordinates (representing matrix entries) for each generator, which we
denote by $g_{ij}$ for 
$i=1,{\ldots} ,n$ and $j=1,2,3,4$. We view $\R (M)$ as a subvariety of
$\C^{4n}$ defined by an ideal $J$ in $\C[g_{11},{\ldots} ,g_{n4}]$.
There are elements $I_\gamma \in \C[g_{11},{\ldots} ,g_{n4}]$ for each
$\gamma \in \pi_1(M)$ such that $I_\gamma (\rho) = \tr \rho (\gamma)$ for
each $\rho \in \R (M)$.

We identify $\m_i$ and $\l_i$ with generators of
$\im(\pi_1(T_i)\to\pi_1(M))$. Thus, $\m_i$ and $\l_i$ are words in the 
generators for $\pi_1(M)$, and we define the
following polynomial equations in the ring
$\C[g_{11},{\ldots},g_{n4},
m_1^{\pm 1},l_1^{\pm 1},{\ldots},m_h^{\pm 1}, l_h^{\pm 1}]$:
\begin{align}
\label{ev:cood m}  I_{\m_i}       &= m_i + m_i^{-1}, \\
\label{ev:cood l}  I_{\l_i}       &= l_i + l_i^{-1}, \\
\label{ev:cood ml} I_{\m_i\l_i} &= m_il_i + m_i^{-1}l_i^{-1},
\end{align}
for $i=1,{\ldots} ,h$. Let $\R_E(M)$ be the variety in
$\C^{4n} \times (\C-\{0\})^{2h}$ defined by $J$ together with the above
equations. For each $\rho \in \R(M)$, the equations
(\ref{ev:cood m}-- \ref{ev:cood ml}) have a solution since commuting
elements of $\SL$ always have a common invariant subspace.
The natural projection $p_1: \R_E(M) \to \R(M)$ is therefore onto, and $p_1$
is a dominating map.

If $(a_1,{\ldots} ,a_{4n}) \in \R (M)$, then there is an action of 
$\Z_2^{h}$ on the resulting points
\begin{equation*}
(a_1,{\ldots} ,a_{4n}, m_1, l_1,{\ldots} ,m_h, l_h) \in \R_E(M)
\end{equation*}
by inverting
both entries of a tuple $(m_i, l_i)$ to $(m_i^{-1}, l_i^{-1})$. 
The group $\Z_2^{h}$ acts transitively on the fibres of $p_1$. The map
$p_1$ is therefore finite--to--one with degree $\le 2^h$. The maximal 
degree is in
particular achieved when the interior of $M$ admits a complete hyperbolic
structure of finite volume, since Theorem \ref{thurston thm} implies that
there are points in $\X_0(M) - \cup_{i=1}^h \{ I^2_{\m_i} = 4\}$.

The \emph{eigenvalue variety} $\Ei (M)$ is the closure of the image of
$\R_E(M)$ under projection onto the coordinates $(m_1,{\ldots} ,l_h)$. It is
therefore defined by an ideal of the ring
$\C [m_1^{\pm 1}, l_1^{\pm 1}, \ldots , m_h^{\pm 1}, l_h^{\pm 1}]$
in $(\C - \{ 0\})^{2h}$.

Note that this construction factors through a variety $\X_E(M)$, which is
the character variety with its coordinate ring appropriately extended, since
we can choose coordinates of $\X (M)$ such
that $I_\gamma$ (as a function on $\X (M)$) is an element of the coordinate
ring of $\X (M)$ for each $\gamma \in \pi_1(M)$.
We let $\te_E : \R_E(M) \to \X_E(M)$ be the natural quotient map,
which is equal to $\te: \R (M) \to \X (M)$ on the first $4n$ coordinates,
and equal to the identity on the remaining $2h$ coordinates.

There also is a restriction map 
$r: \X (M) \to \X(T_1) \times {\ldots} \times \X(T_h)$, which arises
from the inclusion homomorphisms $\pi_1(T_i) \to \pi_1(M)$, and we
therefore denote the map $\X_E(M) \to \Ei (M)$ by $r_E$. Denote the closure
of the image of $r$ by $\X_\partial (M)$. There is the following commuting
diagram of dominating maps:
\begin{center}
$ 
\begin{CD}
      \R_E(M)    @>\te_E >>     \X_E(M)  @>r_E >> \Ei (M)   \\ 
       @Vp_1VV          @Vp_2VV         @Vp_3VV   \\
     \R(M)  @>\te >>     \X(M) @>r >> \X_\partial(M)
\end{CD}
$
\end{center}
Note that the maps $p_1,p_2,t$ and $t_E$ have the property that every point
in the closure of the image has a preimage, and that the maps $p_1, p_2$ and
$p_3$ are all finite--to--one of the same degree.

Recall the construction by Culler and Shalen. We start with a curve
$C \subset \X (M)$ and an irreducible component $R_C$ in $\R (M)$ such that
$\te (R_C) = C$ to obtain the tautological representation 
$\P : \pi_1(M) \to SL_2\C (R_C)$. Let $R'_C$ be an irreducible component of
$\R_E(M)$ with the property that $p_1(R'_C) = R_C$. Then $\C (R'_C)$ is a
finitely generated extension of $\C(R_C)$, and we may think of $\P$ as a
representation  $\P : \pi_1(M) \to SL_2\C (R'_C)$.

If $\overline{r_E\te_E(R'_C)}$ contains a curve $E$, then $\C (R'_C)$ is a
finitely generated extension of $\C (E)$, and hence to each ideal point of
$E$ we can associate essential surfaces using $\P$. Since the eigenvalue of
at least one peripheral element blows up at an ideal point of $E$, the
associated surfaces necessarily have non--empty boundary (see \cite{ccgls}). 

Thus, if there is a closed essential surface associated to an ideal point
of $C$, then either $\overline{r_E\te_E(R'_C)}$ is 0--dimensional, or there
is a neighbourhood $U$ of an ideal point $\xi$ of $R'_C$ such that there 
are (finite) points in $\overline{r_E\te_E(U)} - r_E\te_E(U)$. The later are
called \emph{holes in the eigenvalue variety}, examples of which are given
in \cite{tillus}.


\subsection{Boundary slopes}

We denote the logarithmic limit set of $\Ei (M)$ by $\Ei_\infty (M)$.
There are symmetries of the eigenvalue variety which give rise to symmetries
of its logarithmic limit set. 
If $(m_1, \ldots , m_i, l_i, \ldots , l_h) \in \Ei (M)$,
so is $(m_1, \ldots , m_i^{-1}, l_i^{-1}, \ldots , l_h)$ for any
$i$. Thus, if
$\xi = (x_1,{\ldots} ,x_{2i-1}, x_{2i}, {\ldots} ,x_{2h}) \in\Ei_\infty(M)$,
then $(x_1,{\ldots} ,-x_{2i-1}, -x_{2i}, {\ldots} ,x_{2h}) \in\Ei_\infty(M)$
for $i=1,{\ldots} ,h$. If we factor by these symmetries, we obtain a
quotient of the logarithmic limit set in 
$\RR P^{2h-1}/\Z_2^{h-1} \cong S^{2h-1}$. The quotient map extends to
a map $\Psi : S^{2h -1} \to S^{2h -1}$ of spheres, which has degree $2^h$.
Let 
\begin{equation*}
  T  = \begin{pmatrix} 0 & 1 \\ -1 & 0 \end{pmatrix}
\end{equation*}
and let $T_h$ be the block diagonal matrix with $h$ copies of $T$ on its
diagonal. Then $T_h$ is orthogonal, and its restriction to $S^{2h -1}$ is a
map of degree one. 

\begin{lem}[Boundary slopes] \label{boundary slopes lemma}
  Let $M$ be an orientable, irreducible, compact 3--manifold with non--empty
  boundary consisting of a disjoint union of $h$ tori.
  If $\xi \in \Ei_\infty (M)$ is a point with rational coordinate ratios, then
  there is an essential surface with boundary in $M$ whose projectivised
  boundary curve coordinate is $\Psi (T_h \xi)$.
\end{lem}

Note that the composite $\Psi T_h$ is a smooth map of degree $2^h$.
Since the set of rational points is dense in the logarithmic limit set, the
image $\Psi (T_h \Ei_\infty(M))$ is a closed subset of $\BC (M)$.


\subsection*{Remark}

For a manifold $M$ with boundary consisting of a single torus, Lemma
\ref{boundary slopes lemma} coincides with the result of \cite{ccgls} that
slopes of sides of the Newton polygon of the $A$--polynomial are strongly
detected boundary slopes of $M$. 

In this case, we may normalise $\BC (M)$ and $\Ei_\infty(M)$ such that their
elements are 
pairs of coprime integers. Assume that there is a side of the Newton
polygon of slope $-p/q$ where $p, q \ge 0$ with outward pointing normal vector
$\xi = (q, p)^T$. Then $\xi \in \Ei_\infty(M)$ according to Section
\ref{Spherical Duals}. We now perform matrix multiplication:
\begin{equation*}
  T \xi  = \begin{pmatrix} 0 & 1 \\ -1 & 0 \end{pmatrix}
         \begin{pmatrix} q \\ p \end{pmatrix}
	        = \begin{pmatrix} p \\ -q \end{pmatrix},
		\end{equation*}
and arrive at a point $[p, -q] \in \BC (M)$. The analysis for the remaining
cases is similar, and we have recovered the known relationship.


\subsection*{Proof of Lemma \ref{boundary slopes lemma}}

The coordinate ratios of $\xi$ are rational, and hence Lemma 
\ref{curve finding lemma} provides a curve $C$ in $\Ei (M)$ such
that $\xi \in C_\infty$. From description (\ref{log lim def2}) of elements
in the logarithmic limit set, we know that $\xi$ encodes the definition of a
real--valued valuation on the function field of $C$:
\begin{equation*}
  \xi = (-v(m_1),\ldots ,-v(l_h)).
\end{equation*}
Since the coordinate ratios of $\xi$ are rational, we may rescale the
valuation by $r > 0$, such that for all $i$, $rv(m_i)$ and $rv(l_i)$ are
integers, and such that there is $\alpha \in \Z^{2h} -\{0\}$ with
$rv(m_1^{\alpha_1}\cdots l_h^{\alpha_{2h}}) = 1$. Thus, $rv$ is a normalised,
discrete, rank 1 valuation of $\C (C)$. We finish the proof with standard
arguments. For details see Section 2.1 of \cite{la}, the proof of
Proposition 2.3 of \cite{du2} or Section 5.6 of \cite{ccgls}.
 
The construction by Culler and Shalen yields an action of $\pi_1(M)$ on a
simplicial tree $\tree$ associated to a finitely generated extension $F$ of
$\C (C)$ and to a valuation $w = dv$ for some $d>0$.
The translation length of $\gamma \in \pi_1(M)$ with eigenvalue $e$ (as
an element of $F$) is $\ell (\gamma) = |2 w (e)|$ (see \cite{sh1}, Section
3.9). A dual essential surface $S$ associated to the action can be chosen
such that the geometric intersection number of a loop in $M$ representing
$\gamma$ with $S$ is equal to $\ell (\gamma)$ for all
$\gamma \in \im(\pi_1(T_i) \to \pi_1(M))$.
This implies that a boundary curve of $S$ on $T_i$ is homotopic to one of
$p_i\m_i\pm q_i\l_i$, where $n_ip_i=|2w(l_i)|$ and $n_iq_i=|2w(m_i)|$ and
$p_i,q_i$ are coprime. We need to determine the relative sign.

If $p_i \m_i + q_i \l_i$ represents an element of the subgroup
$\im (\pi_1(\partial S) \to \pi_1(M))$, it is contained in the stabilizer of
an edge, and hence has translation length equal to $0$ since the action is
without inversions.
The elements $\m_i$ and $\m_i^{-1}$ translate points in opposite directions
along their axis. Since $\l_i$ commutes with $\m_i$, they have a common
axis. Thus, the curve $p_i \m_i + q_i \l_i$ can have zero translation length
only if $\m_i$ and $\l_i$ translate in opposite directions. This is the case
if and only if $w(m_i)$ and $w(l_i)$ have opposite signs.
Thus, the relative sign of $p_i$ and $q_i$ is opposite the relative sign of
$w(m_i)$ and $w(l_i)$. We conclude that a boundary curve of $S$ on $T_i$
corresponds to the unoriented homotopy class of the element
$\frac{1}{n_i}(2w(l_i) \m_i - 2w(m_i) \l_i)$.
Thus, the boundary curves of $S$ are encoded by:
\begin{align*}
   & (2w(l_1), - 2w(m_1), {\ldots} , 2w(l_h), - 2w(m_h))\\
  =& 2d(v(l_1), - v(m_1), {\ldots} , v(l_h), - v(m_h))\\
  =& 2d(T_h \xi)
\end{align*}
This completes the proof. \qed


\subsection{Dimension of the eigenvalue and character varieties}

We now deduce some bounds on the dimension of the eigenvalue and the
character variety. 

\begin{pro} \label{ev:dim pro}
Let $M$ be an orientable, irreducible, compact  3--manifold with boundary
consisting of $h$ tori. Then $\dim_\sC \Ei (M) \le h$. 
  
Moreover, if the interior of $M$ admits a complete hyperbolic structure of
finite volume, then $\dim_\sC \Ei (M) = h$. 
\end{pro}  
\begin{proof}
We have $h>\dim_\RR \BC (M) \ge \dim_\RR \Ei_\infty = \dim_\sC \Ei (M) -1$,
where the first inequality is due to Theorem \ref{ev:dim bc}, the second to
Lemma  \ref{boundary slopes lemma}, and the equality to Theorem
\ref{Bergman-Bieri-Groves}. It follows that $h \ge \dim_\sC \Ei (M)$.

To prove the second statement, it is enough to show that under the
assumptions, there is a subvariety of $\Ei (M)$ of dimension $\ge h$.
Theorem \ref{thurston thm} asserts that there is a component $\X_0(M)$ of
$\X (M)$ of complex dimension equal to $h$. Thurston shows in \cite{t},
Section 5.8, that the map $\X_0(M) \to \C^h$ defined by
$\chi_\rho \to (I_{\m_1}(\rho),{\ldots} ,I_{\m_h}(\rho))$ maps a small open
neighbourhood of $\chi_0$ to an open neighbourhood of
$(\pm 2,{\ldots} ,\pm 2)$. This is only possible if
$I_{\m_i} = \tr\rho(\m_i)$ are algebraically independent over $\C$ as
elements of $\C [X_0]$.
There is a component $X'_0$ of $\X_E(M)$ with the
same properties, since $p_2$ is a finite dominating map and the induced
homomorphism $p_{2*}: \C [X_0] \to \C [X'_0]$ is injective.

Let $E_0 = \overline{r_E\te_E(X'_0)}$. If the functions
$m_i+m_i^{-1}$, $i=1,{\ldots} ,h$, are algebraically dependent over $\C$ as
elements of $\C[E_0]$, then they are also algebraically dependent over $\C$
as elements of $\C [X'_0]$ since we have a ring homomorphism
$\C[E_0] \to \C [X'_0]$. But the identities (\ref{ev:cood m}) imply that then
the functions $I_{\m_i}$ are algebraically dependent over $\C$ as elements
of $\C [X'_0]$, contradicting the fact that $p_{2*}$ is injective.

Thus, the functions $m_i+ m_i^{-1}$, $i=1,{\ldots} ,h$, are algebraically
independent over $\C$ as elements of $\C[E_0]$, and hence
$\dim_\sC E_0 \ge h$ according to (\cite{clo}, 9.5 Theorem 2). This
completes the proof.
\end{proof}

\begin{cor}
Let $M$ be an orientable, irreducible, compact 3--manifold with boundary
consisting of $h$ tori. 

If $M$ contains no closed essential surface, then
$\dim_\sC \X (M) = \dim_\sC \Ei (M) \le h$.
\end{cor}
\begin{proof}
We have $\dim_\sC \X (M) = \dim_\sC \X_E(M)$. Assume there is a component
$X$ of $\X_E(M)$ such that $\dim X > \dim \overline{r_E(X)}$. Then
\cite{mu}, Theorem (3.13), implies that there is a subvariety $V$ of
dimension $\ge 1$ in $X$ which maps to a single point of $\Ei (M)$.
Associated to an ideal point of a curve in $V$, we find a closed essential
surface since the traces of all peripheral elements are constant on $V$.
This contradicts our assumption on $M$. Thus, for all components $X$ of 
$\X (M)$, we have  $\dim X = \dim \overline{r_E(X)}$, and this proves the
claim. 
\end{proof}


\subsection{Remark}

Let $M$ be an orientable, irreducible, compact 3--manifold with boundary
consisting of $h$ tori, such that the interior of $M$ admits a complete
hyperbolic structure of finite volume. The proof of Proposition \ref{ev:dim
pro} implies that the restriction of $r$ to the Dehn surgery component,
$r_0: X_0 \to \X_\partial (M)$, is finite--to--one onto its image. The degree
of $r_0$ depends upon $H^1(M, \Z_2)$, so let us consider the corresponding map
$\overline{r}_0: \overline{X}_0 \to \PX_\partial (M)$. In the case where
$h=1$, Dunfield has shown in \cite{du2} that $\overline{r}_0$ has degree 1
onto its image using Thurston's Hyperbolic Dehn Surgery Theorem and a Volume
Rigidity Theorem (Gromov--Thurston--Goldman). It would be interesting to
generalise Dunfield's proof to arbitrary $h$, showing that $\overline{r}_0$
is always a birational isomorphism onto its image.


\subsection{$\mathbf{\PSL}$--eigenvalue variety} 

Analogous to $\Ei(M)$, one can define a $\PSL$--eigenvalue variety
$\PEi(M)$, since the function $I^2_\gamma : \PX (M) \to \C$ defined by
$I^2_\gamma (\prho) = (\tr\prho(\gamma))^2$ is regular (i.e. polynomial) for
all $\gamma\in\pi_1(M)$ (see \cite{bozh}). Thus, we construct a variety
$\PR_E(M)$ using the relations
\begin{align*}
I^2_{\m_i} &= M_i + 2 + M_i^{-1}, \\
I^2_{\l_i} &= L_i + 2+L_i^{-1}, \\
I^2_{\m_i\l_i} &= M_iL_i + 2+M_i^{-1}L_i^{-1},
\end{align*}
for $i=1,{\ldots} ,h$. 

Consider the representation of $\Z \oplus \Z$ into $\PSL$ generated by the
images of
\begin{equation} \label{four group}
    \begin{pmatrix} i&0 \\ 0&-i \end{pmatrix}
    \quad \text{and} \quad
    \begin{pmatrix} 0&1 \\ -1&0 \end{pmatrix}.
\end{equation}
In $\PSL$, the image of this representation is isomorphic to
$\Z_2 \oplus \Z_2$, but the image of any lift to $\SL$
is isomorphic to the quaternion group $Q_8$ (in Quaternion group notation).
If $\prho \in \PR (M)$ restricts to such an irreducible abelian
representation on a boundary torus $T_i$, then the equations
$M_i + 2 + M_i^{-1}=0$, $L_i + 2+L_i^{-1}=0$,
$M_iL_i + 2+M_i^{-1}L_i^{-1}=0$ have no solution. In particular, the
projection $\overline{p}: \PR_E(M) \to \PR(M)$ may
not be onto.

The closure of the image of $\PR_E(M)$ onto the coordinates
$(M_1,L_1,{\ldots},M_h ,L_h)$ is denoted by $\PEi (M)$ and called the
\emph{$\PSL$-eigenvalue variety}. The relationship between ideal points of
$\PEi (M)$ and strongly detected boundary curves is the same as for the
$\SL$--version, since the proof of Lemma \ref{boundary slopes lemma} applies
without change:

\begin{lem} \label{psl boundary slopes lemma}
  Let $M$ be an orientable, irreducible, compact 3--manifold with non--empty
  boundary consisting of a disjoint union of $h$ tori.
  If $\xi \in \PEi_\infty (M)$ is a point with rational coordinate ratios, then
  there is an essential surface with boundary in $M$ whose projectivised
  boundary curve coordinate is $\Psi (T_h \xi)$.
\end{lem}


\section{Torus knots}

In general, one needs a computer to calculate the $A$--polynomial. It is shown
in \cite{ccgls} that if $\k$ is a nontrivial $(p,q)$--torus knot, then
$A_\k(l,m)$ is divisible by $lm^{pq} + 1$ (where we follow the convention
that if $M$ is the complement of a knot $\k$ in $S^3$, then $\{\m, \l\}$ is
a standard peripheral system).

\begin{pro}\footnote{I thank Walter Neumann for pointing out that an earlier
argument for $(\alpha, 1)$--two bridge knots could be generalised to torus
knots.} \label{torus knot A-polys} 
  Let $\k$ be a $(p,q)$--torus knot. If $p=2$ or $q=2$, then
  $A_\k(l,m) = (l-1)(lm^{pq} + 1)$, otherwise
  $A_\k(l,m) = (l-1)(lm^{pq} + 1)(lm^{pq}-1)$.
\end{pro}

\begin{proof}
  The fundamental group of $S^3 - \k$
  is presented by $\Gamma = \langle u,v \mid u^p = v^q\rangle$, and
  a standard peripheral system is given by
  $\m = u^nv^m$ and $\l = u^p\m^{-pq}$, where $mp+nq = 1$. These facts
  can be found in \cite{bz} on page 45. The factor $l-1$ arises from
  reducible representations, and all factors arising from
  components containing irreducible representations are to be determined.

  The element $u^p$ is in the centre of $\Gamma$, since it is identical to
  $v^q$ and hence commutes with both generators.
  Thus the image of $u^p$ is in the
  center of $\rho (\Gamma )$. Since two commuting elements have a common
  eigenvector, $\rho (u)^p $ has a common eigenvector with $\rho
  (u)$ and $\rho (v)$ respectively. If the representation is
  irreducible, these eigenvectors have to be distinct. Thus, after
  conjugation it may be assumed that the generators map to upper and
  lower triangular matrices respectively, and the central element is
  represented by a diagonal matrix. Direct matrix computations show that the
  commutativity with either of the generators gives $\rho (u)^p = \pm E$. 

  Assume that $\rho (u)^p = -E$. The relation $\rho (\l ) = -\rho
  (\m )^{-pq}$, which is equivalent to $\rho (\l \m^{pq})= -E$,
  then implies the equation $lm^{pq} = -1$. This is the curve obtained in
  \cite{ccgls} by sending $u$ and $v$ to noncommuting elements of $SL_2(\C )$
  of order exactly $2p$ and $2q$ respectively.

  If the image of $\rho (u)^p$ is trivial, then $\rho (v)^q = \rho (u)^p = E$.
  If $p$ or $q$ equals $2$, this implies that the image of one of the
  generators is $\pm E$. But this yields that $\rho (\Gamma )$ is
  abelian, and therefore contradicts the irreducibility assumption. This
  completes the proof of the first assertion.

  Now assume that neither of $p$ and $q$ is equal to two.
  As above, curves of irreducible representations are obtained by sending
  $u$ and $v$ to noncommuting elements of $SL_2(\C )$ of order exactly $p$
  and $q$ respectively. Any of these curves yields a
  component of the eigenvalue variety defined by the equation $lm^{pq} = 1$,
  and this
  finishes the proof of the proposition. In fact, with a more precise
  analysis one could count the number of 1--dimensional curves in the
  character variety of a torus knot. 
\end{proof}

Similarly to the above, we can compute a generator for the defining
(principal) ideal of the $\PSL$--eigenvalue variety of a torus knot. We call
this generator the $\overline{A}$--polynomial, and its variables have the
unfortunate names $L$ and $M$. Thus, $\overline{A}_M(L,M)$ is the
$\overline{A}$--polynomial associated to the manifold $M$.

\begin{pro} \label{torus knot A_bar-polys}
  Let $\k$ be a $(p,q)$--torus knot. Then
  $\overline{A}_\k(L,M) = (L-1)(LM^{pq} - 1)$.
\end{pro}

\begin{proof}
It follows from the lemma in Section 6 of \cite{ccgls} that all
$\PSL$--representations of a knot group lift to $\SL$. Thus, all factors of
the $A$--polynomial arise from factors of the $\overline{A}$--polynomial and
vice versa. 

The eigenvalue of the meridian appears with even powers in the (unfactorised)
$A$--polynomial of $\k$, and the sign of the longitude's eigenvalue is
uniquely  determined since $\l$ is null--homologous. 
Thus, $A_\k(l,m) = A_\k(l,-m)$, but $A_\k(l,m) = 0$ does not imply 
$A_\k(-l,m)=0$ in
general. However, since 
$\overline{A}_\k(l^2, m^2) = \overline{A}_\k((-l)^2, m^2)$, it follows that
the variety defined by $\overline{A}_\k(l^2,m^2)=0$ is the locus of
$A_\k(l,m)A_\k(-l,m)=0$.
Expanding $A_\k(l,m)A_\k(-l,m)$, substituting $l^2 = L$ and $m^2 = M$ and
deleting repeated factors therefore gives the result as stated.
\end{proof}


\subsection*{Acknowledgements}

The contents of this paper forms part of my Ph.D. thesis, and I thank my
supervisors Craig Hodgson and Walter Neumann for their constant support
and inspiration. This research was supported in part by an International
Postgraduate Research Scholarship by the Commonwealth of Australia
Department of Education, Science and Training. 


{\small


{\small Department of Mathematics and Statistics,
The University of Melbourne,
VIC 3010, Australia
}
}

\end{document}